\documentclass[a4paper, english,11pt, reqno]{amsart}

\usepackage{graphicx} 
\usepackage{amsmath}
\usepackage{amssymb}
\usepackage{comment}

\numberwithin{equation}{section}
\numberwithin{figure}{section}

\newtheorem{thm}{\protect\theoremname}
 
  \newtheorem{lemma}[thm]{\protect\lemmaname}
    \newtheorem{prop}[thm]{\protect\propname}

\newtheorem{df}[thm]{Definition}

\numberwithin{thm}{section}

\newtheorem{cor}[thm]{Corollary}

\newtheorem{theorem}{Theorem}

\newtheorem{rem}{Remark}

\makeatother

\usepackage{babel}
\providecommand{\propname}{Proposition}

\providecommand{\lemmaname}{Lemma}
\providecommand{\theoremname}{Theorem}

\newcommand{\MA}{\operatorname{MA}}
\newcommand{\E}{\operatorname{E}}
\newcommand{\NA}{\operatorname{NA}}
\newcommand{\di}{\operatorname{div}}
\newcommand{\PSH}{\operatorname{PSH}}
\newcommand{\com}{\operatorname{com}}

\title[Two Calabi-Yau theorems for degenerations]{Two Calabi-Yau theorems for degenerations of compact Kähler manifolds}
\author{David Witt Nyström}
\date{September 2025}

\begin{document}

\maketitle

\begin{abstract}
We discuss two closely related Calabi-Yau theorems for degenerations of compact Kähler manifolds. The first is a Calabi-Yau theorem for big test configurations, that generalizes a result in \cite{WN24}. It follows from recent joint work with Mesquita-Piccione \cite{MW25}, but is here given a more direct proof. The second result is a Calabi-Yau theorem for a wider class of degenerations, formulated in the language of non-Archimedean Kähler geometry. It was first proved in the algebraic setting by Boucksom-Jonsson \cite{BJ22}, building on earlier work of Boucksom-Favre-Jonsson \cite{BFJ15}, while the general Kähler case was established in \cite{MW25}. Our main focus here is on the connection between these results and the theory of big cohomology classes and their volumes.    
\end{abstract}

\section{Introduction}
Let $(X,\omega)$ be a compact Kähler manifold and let $\alpha:=\{\omega\}\in H^{1,1}(X,\mathbb{R})$ be the associated Kähler class. To understand $(X,\omega)$ or $(X,\alpha)$, it is often helpful to study their degenerations. As an example, by the work of Chen-Cheng \cite{CC21a,CC21b}, we know that the existence of a constant scalar curvature Kähler (cscK) metric in $\alpha$ can be detected using geodesic rays, which are degenerations of $(X,\omega)$ of a special kind. There is also a strong link between geodesic rays and certain degenerations of $(X,\alpha)$ called test configurations. Test configurations appear in the formulation of Yau-Tian-Donaldson (YTD) conjecture, which aims to give a numerical criterion for the existence of canonical metrics, such as cscK metrics, in a given Kähler class. After e.g. important contributions by Tian, see e.g. \cite{Tia97,Tia15}, the YTD conjecture was famously proved by Chen-Donaldson-Sun in the case when $X$ is Fano and $\alpha=c_1(X)$ \cite{CDS15a,CDS15b,CDS15c}. Recently, Boucksom-Jonsson \cite{BJ25} proved a version of the YTD conjecture in the algebraic setting, building on earlier work of Chi Li \cite{Li22,Li23a}, and in \cite{MW25} Mesquita-Piccione and I extended some of latter results to the Kähler setting. At the same time, Darvas-Zhang \cite{DZ25} proved a different YTD correspondence, that is also valid in the Kähler setting. 

The aim of this paper is to discuss two closely related Calabi-Yau theorems for degenerations of $(X,\alpha)$, Theorem \ref{thm:A} and Theorem \ref{thm:main}. Theorem \ref{thm:A} is a Calabi-Yau theorem for so-called big test configurations, and it generalizes a result in \cite{WN24}. Theorem \ref{thm:main} is a Calabi-Yau theorem for a wider class of degenerations, formulated using the language of non-Archimedean Kähler geometry. Theorem \ref{thm:main} plays a key role in the variational approach to the YTD conjecture, see e.g. \cite{BBJ21,Li22}, and was first proved in the algebraic case (i.e. when $X$ is projective and $\alpha=c_1(L)$ for some ample $(\mathbb{R}-$)line bundle $L$) by Boucksom-Jonsson \cite{BJ22}, building on earlier work of Boucksom-Jonsson-Favre  \cite{BFJ15}. The general Kähler case of Theorem \ref{thm:main} was recently proved by Mesquita-Piccione and myself \cite{MW25}. 

\subsection*{Organization}
The paper is organized as follows.

We first review the classical Calabi-Yau theorem in Section \ref{sec:CYT}. 

In Section \ref{sec:tc} we discuss Kähler and big test configurations and give a geometric interpretation of their associated Monge-Amp\`ere measures. Here we also formulate Theorem \ref{thm:A}. 

In Section \ref{sec:nAKg} we introduce the notions of non-Archimedean Kähler geometry needed to formulate Theorem \ref{thm:main}. 

In Section \ref{sec:big} we give a brief introduction to the theory of big cohomology classes. In particular, we discuss the conjectural duality between the pseudoeffective and the movable cone, the related transcendental Morse inequality, and how this connects to differentiability properties of the volume function. Finally, we state a result from \cite{WN24} about restricted volumes, here called Corollary \ref{cor:WN}, which will be crucial for what follows.

In Section \ref{sec:direct} we use the general results of Section \ref{sec:big}, and in particular Corollary \ref{cor:WN}, to give a direct proof of Theorem \ref{thm:A}.

In Section \ref{sec:proof} we discuss the proof of Theorem \ref{thm:main}, which uses the variational method described in \cite{BFJ15,BJ22,BJ23a} for the algebraic case. A key step is to establish the so-called orthogonality property, see Definition \ref{def:orth}, and we show how this can be accomplished using the results from Section \ref{sec:direct}.

\subsection*{Acknowledgements}
My understanding of this topic has been shaped by interactions with many great people, such as Tristan Collins, Tamás Darvas, Ruadhaí Dervan, Eleonora Di Nezza, Vincent Guedj, Jakob Hultgren, Hoang-Chinh Lu, Mihai Păun, Rémi Reboulet, Yanir Rubinstein, Valentino Tosatti, Antonio Trusiani, Mingchen Xia, and in particular Sébastien Boucksom and Mattias Jonsson, who developed a lot of the mathematics that is discussed in this paper. Jean-Pierre Demailly, who sadly is not with us anymore, was another great influence. I also want to express my deep-felt gratitude to my former advisors Robert Berman and Bo Berndtsson, my long-time collaborator Julius Ross, and my more recent Calabi-Yau collaborator Pietro Mesquita-Piccione. 

\section{The classical Calabi-Yau theorem} \label{sec:CYT}

Here and in the rest of this paper, $(X,\omega)$ will be a compact Kähler manifold of complex dimension $n$, with the associated Kähler class denoted by $\alpha:=\{\omega\}\in H^{1,1}(X,\mathbb{R})$. We let $V:=\int_X \omega^n=\alpha^n$ denote the volume.

One of the most central results in Kähler geometry is the original Calabi-Yau theorem, due to Yau \cite{Yau77,Yau78}.

\begin{thm} \label{thm:CY1}
For any volume form $dV$ on $X$ such that $\int_XdV=V$, there is a unique K\"ahler form $\omega'\in\alpha$ such that $(\omega')^n=dV$.
\end{thm}

Note that by the $dd^c$-lemma, any Kähler form $\omega'\in \alpha$ can be written as $\omega'=\omega+dd^c\phi$ for some smooth function $\phi$. We say that a smooth function $\phi$ is a Kähler potential (with respect to $\omega$) if $\omega_{\phi}:=\omega+dd^c\phi$ is Kähler. The set of Kähler potentials is denoted by $\mathcal{H}_{\omega}$, and we also let $\mathcal{H}_{\omega,0}$ denote the set of Kähler potentials whose supremum is zero. The Monge-Amp\`ere measure of $\phi$ is defined as $\MA_{\omega}(\phi):=V^{-1}\omega_{\phi}^n$, thus Theorem \ref{thm:CY1} can be restated as saying that for any volume form $dV$ such that $\int_X dV=1$, the Monge-Amp\`ere equation 
\begin{eqnarray*}
\MA_{\omega}(\phi)=dV
\end{eqnarray*}
has a unique solution in $\mathcal{H}_{\omega,0}$. 

There are also variations of Theorem \ref{thm:CY1}, where the domain of definition of the Monge-Amp\`ere operator is enlarged, see e.g. \cite{Koł98} and \cite[Theorem A]{BBGZ13}. To formulate the second of these variations, we need to recall some basic concepts of pluripotential theory.

A decreasing limit $\psi$ of K\"ahler potentials (with respect to $\omega$), not identically equal to $-\infty$, is said to be $\omega$-psh, and the set of $\omega$-psh functions is denoted by $\PSH_{\omega}$. 

The energy of a K\"ahler potential is defined as $$\E_{\omega}(\phi):=\frac{V^{-1}}{n+1}\sum_{j=0}^n\int_X \phi(\omega+dd^c\phi)^j\wedge \omega^{n-j}.$$ The energy of $\psi\in \PSH_{\omega}$ is defined as the infimum of the energy of all $\phi\in \mathcal{H}_{\omega}$ such that $\phi\geq \psi$, and the space of finite energy potentials is defined as $$\mathcal{E}^1_{\omega}:=\{\psi\in \PSH_{\omega}: \E_{\omega}(\psi)>-\infty\}.$$ We also let $\mathcal{E}^1_{\omega,0}:=\{\phi\in \mathcal{E}^1_{\omega} : \sup \phi=0\}$.

The Monge-Amp\`ere measure of a K\"ahler potential $\phi$ is defined as $$\MA_{\omega}(\phi):=V^{-1}(\omega+dd^c\phi)^n,$$ and there is a natural extension of the Monge-Amp\`ere operator to $\mathcal{E}^1_{\omega}$. 

The (dual) energy $\E^{\vee}_{\omega}(\mu)$ of a Radon probabiltiy measure is defined as $$\E^{\vee}_{\omega}(\mu):=\sup\left\{\E_{\omega}(\phi)-\int_X\phi\, \mathrm d\mu: \phi\in \mathcal{E}^1_{\omega}\right\},$$ those with $\E^{\vee}_{\omega}(\mu)<\infty$ giving us the space of finite energy measures $\mathcal{M}^1_{\omega}$.

In \cite{BBGZ13} Berman-Boucksom-Guedj-Zeriahi used variational methods to prove the following Calabi-Yau theorem for finite energy potentials. 

\begin{thm} \label{thm:CY}
The Monge-Amp\`ere operator is a bijection between $\mathcal{E}^1_{\omega,0}$ and $\mathcal{M}^1_{\omega}$.
\end{thm}

\section{A Calabi-Yau theorem for big test configurations} \label{sec:tc}

\subsection{Test configurations and models}

Let us first recall the definition of a test configuration, which in the general Kähler setting goes back to \cite{DR17} and \cite{SD18}.

\begin{df} A (smooth dominating) test configuration $(\mathcal{X},A+D)$ of $(X,\alpha)$ consists of the following data:
\begin{enumerate}
\item a compact K\"ahler manifold $\mathcal{X}$ together with a surjective map $\pi:\mathcal{X}\to X\times \mathbb{P}^1$ such that $\pi: \mathcal{X}\setminus |\mathcal{X}_0|\to X\times (\mathbb{P}^1\setminus \{0\})$ is a biholomorphism, where $\mathcal{X}_0:=\pi_{\mathbb{P}^1}^*([0])$ is the zero divisor and $|\mathcal{X}_0|:=\pi_{\mathbb{P}^1}^{-1}(\{0\})$ is the zero fiber,
\item a lift of the standard $\mathbb{C}^*$-action on $X\times \mathbb{P}^1$ to $\mathcal{X}$ making $\pi$ equivariant,
\item a class $A+D\in H^{1,1}(\mathcal{X},\mathbb{R})$ where $A:=\pi_X^*(\alpha)$ and $D$ is a vertical divisor, i.e. a divisor supported on $|\mathcal{X}_0|$ (for convenience we do not distinguish between $D$ and its cohomology class).
\end{enumerate}

\end{df}

We also call $\mathcal{X}$ a model of $X$. If $\mathcal{X}$ is a model and $Y\subseteq |\mathcal{X}_0|$ is a $\mathbb{C}^*$-invariant submanifold, the blow-up of $Y$ in $\mathcal{X}$ is a new model. A simple way to produce models is thus to start from the trivial model $X\times \mathbb{P}^1$, and then to iterate this blow-up procedure.

We say that the model or test configuration is SNC if $\mathcal{X}_0^{red}$ has simple normal crossings, i.e. if $\mathcal{X}_0=\sum_i b_iE_i$ where the $E_i$:s are smooth hypersurfaces that intersect transversely. From now on all models and test configurations are assumed to be SNC, unless specifically stated otherwise. 

$\sum_{i=1}^Nb_iE_i$ will always denote the decomposition of $\mathcal{X}_0$ into its weighted irreducible components, and we also let $\mathcal{X}_0^{\com}:=\{E_i : 1\le i\le N\}$.

\subsection{Kähler test configurations}

We say that the test configuration $(\mathcal{X},A+D)$ is K\"ahler if $A+D$ is Kähler, meaning that $A+D$ contains a Kähler form. 

Let thus $(\mathcal{X},A+D)$ be a Kähler test configuration and $\Omega$ a Kähler form in $A+D$. 

Let also $X_{\tau}\cong X$ denote the fiber of over $\tau\in \mathbb{P}^1\setminus\{0\}$, and let $\omega_{\tau}:=\Omega_{|X_{\tau}}$. Since $D$ does not intersect $X_{\tau}$ we have that $\omega_{\tau}\in A_{|X_{\tau}}=\alpha$, and thus $(X_{\tau},\omega_{\tau})_{\tau\in \mathbb{P}^1\setminus\{0\}}$ is a family of Kähler manifolds, all with volume $V$. But as $\tau\to 0$, $(X_{\tau},\omega_{\tau})$ degenerates to $(\mathcal{X}_0, \Omega_{|\mathcal{X}_0})$, and the proportion of the fixed volume $V$ that goes into the weighted component $b_iE_i$ is given by $$V^{-1}b_i\int_{E_i}\Omega^n=V^{-1}(A+D)^n\cdot (b_iE_i).$$ 

This motivates the following definition.

\begin{df}
    The Monge-Amp\`ere measure of a Kähler test configuration $(\mathcal{X},A+D)$ is the measure on $\mathcal{X}_0^{\com}$ defined as $$\MA(\mathcal{X},A+D):=V^{-1}\sum_{i=1}^N(A+D)^n\cdot (b_iE_i)\delta_{E_i}.$$
\end{df}

Thus the Monge-Amp\`ere measure encodes how the volume is distributed among the weighted components of $\mathcal{X}_0$. As we will see in Section \ref{sec:nAKg}, this definition of $\MA(\mathcal{X},A+D)$ comes from non-Archimedean Kähler geometry, where it plays a central role.

Since 
\begin{eqnarray*}
\sum_{i=1}^N(A+D)^n\cdot (b_iE_i)=(A+D)^n\cdot \mathcal{X}_0=\\=(A+D)^n\cdot X_1=A^n\cdot X_1=\alpha^n=V,
\end{eqnarray*}
we see that $\MA(\mathcal{X},A+D)$ is a probability measure on $\mathcal{X}_0^{\com}$, which gives each point a non-zero mass. If $\mu$ is such a measure on $\mathcal{X}_0^{\com}$, i.e. if $\mu=\sum_{i=1}^n a_i\delta_{E_i}$ with $a_i>0$ and $\sum_i a_i=1$, can we always find a Kähler test configuration $(\mathcal{X},A+D)$ which solves the Monge-Amp\`ere equation $$\MA(\mathcal{X},A+D)=\mu?$$ 

The answer is no. A very simple example where it fails was given in \cite{WN24}, and we will quckly recall that here. Let $X:=\mathbb{P}^1\times \mathbb{P}^1$, $\alpha:=\{\pi^*_1\omega_{FS}+\pi^*_2\omega_{FS}\}$ and $\mathcal{X}:=Bl_{(0,0)}X$. Then $\mathcal{X}_0=E_1+E_2$ where $E_1$ is the proper transform of $\mathbb{P}^1\times\{0\}$ and $E_2$ is the exceptional divisor of the blow-up. In this case, as shown in \cite{WN24}, we can find a Kähler test configuration $(\mathcal{X},A+D)$ such that $\MA(\mathcal{X},A+D)=a_1\delta_{E_1}+a_2\delta_{E_2}$ if and only if $1/2<a_1<1, 0<a_2<1/2$ and $a_1+a_2=1$. 

\subsection{Big test configurations}

To be able to solve the Monge-Amp\`ere equation for arbitrary probability measures, we are thus forced to consider a larger class of test configurations. 

We say that a test configuration $(\mathcal{X},A+D)$ is big if the class $A+D$ is big, i.e. if $A+D$ can be written as the sum of a Kähler class and a pseudoeffective class (for more details on big classes see Section \ref{sec:big}). Note that in the algebraic setting with $A+D=c_1(\mathcal{L})$, $A+D$ is big precisely when $\mathcal{L}$ is a big line bundle, i.e. when $h^0(\mathcal{X},\mathcal{L}^k)$ grows like a positive multiple of $k^{n+1}$.

\begin{df}
The Monge-Amp\`ere measure of a big test configuration $(\mathcal{X},A+D)$ is the measure on $\mathcal{X}_0^{\com}$ defined as $$\MA(\mathcal{X},A+D):=V^{-1}\sum_{i=1}^Nb_i\langle(A+D)^n\rangle_{\mathcal{X}|E_i}\delta_{E_i}.$$
\end{df}

Here $\langle(A+D)^n\rangle_{\mathcal{X}|E_i}$ denotes the restricted volume of $A+D$ along $E_i$, which is equal to $(A+D)^n\cdot E_i$ when $A+D$ is Kähler (see Section \ref{sec:big} for the definition). In contrast to the intersection number $(A+D)^n\cdot E_i$, the restricted volume is always nonnegative, which in particular means that $\MA(\mathcal{X},A+D)$ is a positive measure.  

Let us now give a geometric interpretation of $\MA(\mathcal{X},A+D)$.

A big class $A+D$ on $\mathcal{X}$ that is not Kähler will obviously not contain any Kähler forms. However, it will contain many closed positive currents $\Omega$ that are smooth Kähler forms away from the so-called non-Kähler locus $E_{nK}:=E_{nK}(A+D)$ of $A+D$ (see Section \ref{sec:big} for the definition), which in this case will be an analytic subset of $|\mathcal{X}_0|$. Given such an $\Omega$ we let $\omega_{\tau}:=\Omega_{|X_{\tau}}$ and thus $(X_{\tau},\omega_{\tau})_{\tau\in \mathbb{P}^1\setminus\{0\}}$ is a family of Kähler manifolds, all with volume $V$. As $\tau\to 0$, $(X_{\tau},\omega_{\tau})$ degenerates to $(\mathcal{X}_0, \Omega_{|\mathcal{X}_0})$. Since $\Omega$ can have  singularities it will now typically happen that $$\sum_ib_i\int_{E_i\setminus E_{nK}}\Omega^n<V,$$ which tells us that as $\tau\to 0$, some of the volume of $(X_{\tau},\omega_{\tau})$ concentrates along the non-Kähler locus $E_{nK}$. Now $\langle(A+D)^n\rangle_{\mathcal{X}|E_i}$ is equal to the supremum of all possible integrals $\int_{E_i\setminus E_{nK}}\Omega^n$ where $\Omega$ is chosen as above, and by monotonicity one can find a sequence of such currents $\Omega_k$ so that for all $i$, $\int_{E_i\setminus E_{nK}}\Omega_k^n$ will converge to $\langle(A+D)^n\rangle_{\mathcal{X}|E_i}$.

In Section \ref{sec:direct} we will prove the following.

\begin{thm} \label{thm:prob}
The Monge-Amp\`ere measure of a big test configuration is always a probability measure.   
\end{thm}

Note that because of the singularities of the currents $\Omega$, this does not follow directly from the fact that $\{X_1\}=\{\mathcal{X}_0\}$, as it did in the Kähler case. Indeed, it is a highly nontrivial result, closely related to the so-called transcendental Morse inequality (see Sections \ref{sec:big} and \ref{sec:direct}). 

We thus see that the Monge-Amp\`ere measure $\MA(\mathcal{X},A+D)$ encodes how the volume of $(X_{\tau},\omega_{\tau})$ distributes among the weighted components as $\tau\to 0$, at least up to discrepancies due to the singularities of $\Omega$, but as a consequence of Theorem \ref{thm:prob}, these discrepancies can be made arbitrarily small. 

The Calabi-Yau theorem for big test configurations now says that the associated Monge-Amp\`ere equation always is solvable.

\begin{theorem} \label{thm:A}
If $\mathcal{X}$ is a model, then for any probability measure $\mu$ on $\mathcal{X}_0^{\com}$ one can find a big test configuration $(\mathcal{X},A+D)$ such that $$\MA(\mathcal{X},A+D)=\mu.$$
\end{theorem}

In other words, using big test configurations it is possible to prescribe how the volume distributes among the weighted components of the zero divisor.

In the algebraic setting Theorem A follows e.g. from \cite[Corollary 7.16]{BJ25a}. In the general Kähler setting it is a consequence of \cite[Proposition 8.3.2]{MW25}, which is proved by invoking the more general Theorem B. However, there is also a more direct proof. Namely, in \cite{WN24} I proved the special case of Theorem \ref{thm:A} when $\mu$ is a Dirac measure, and as we will see in Section \ref{sec:direct}, the full statement can be proved using a refinement of the argument in \cite{WN24}.

\section{A non-Archimedean Calabi-Yau theorem} \label{sec:nAKg}
As already noted, the definition of the Monge-Amp\`ere measure of an ample/Kähler test configuration originated in non-Archimedean Kähler geometry, or non-Archimedean pluripotential theory which it is also called.

This theory, where common notions of Kähler geometry are given analogues in a non-Archimedean context, was initiated by Kontsevich-Tscinkel \cite{KT01}, and has since been extensively developed by Boucksom, Chambert-Loir, Ducros, Favre, Jonsson and others (see e.g. \cite{BFJ15,BJ22,BJ23a,CD12} and references therein). Notably, the theory underlies the variational approach to the Yau-Tian-Donaldson conjecture, and is also the basis for Yang Li's approach to the Strominger-Yau-Zaslow conjecture, see e.g. \cite{YaLi22,YaLi23,HJMM24}.

Originally, non-Archimedean Kähler geometry was formulated only in the algebraic setting. Recently though, in the special case relevant to the Yau-Tian-Donaldson conjecture (i.e. when the non-Archimedean field is given by $\mathbb{C}$ with its trivial norm), Darvas-Xia-Zhang \cite{DXZ23} and Mesquita-Piccione \cite{MP24} proposed two somewhat different ways of formulating the theory for compact Kähler manifolds. Here we will use the version presented in \cite{MP24}, which indeed is close to the original formulation in \cite{KT01}. In addition to \cite{MP24} and \cite{MW25}, key references include \cite{BFJ15,BJ22,BJ23a}.

\subsection{Tropical analytification} Let $\mathcal{I}_X$ denote the set of coherent ideal sheaves on $X$. A non-constant function $v: \mathcal{I}_X\to [0,\infty]$ is called a semivaluation if for any $I,J\in\mathcal{I}_X$ we have that $v(IJ)=v(I)+v(J)$ and $v(I+J)=\min(v(I),v(J))$.
The tropical analytification $X^{\NA}$ of $X$ is defined as the set of semivaluations given the topology of pointwise convergence. $X^{\NA}$ then  becomes a compact Hausdorff space.

\subsection{Divisorial points} Let now $\mathcal{X}$ be a model (recall that all our models are assumed to be SNC) with zero divisor $\mathcal{X}_0=\sum_i b_iE_i$, as defined in Section \ref{sec:tc}. Then for each irreducible component $E_i$ there is an associated semivaluation $$v_{E_i}(I):=\min\{b_i^{-1}\textrm{ord}_{E_i}(f\circ \pi_X): f\in I(U), U\subseteq X\}.$$
Such semivaluations are called \emph{divisorial  valuations} and the set of divisorial valuations/points in $X^{\NA}$ is denoted by $X^{\di}$. Thus we can think of $\mathcal{X}_0^{\com}$ as a subset of $X^{\di}$. Note that if  $\mu: \mathcal{X}'\to \mathcal{X}$ is a dominating model and $E_i'$ is the proper transform of $E_i$ then $v_{E_i}=v_{E_i'}$. A crucial fact is that $X^{\di}$ is dense in $X^{\NA}$. 

\subsection{Dual complexes} Let $\mathcal{X}$ be a model with $\mathcal{X}_0=\sum_{i\in I} b_iE_i$. To each subset $J\subseteq I$ and irreducible component $Z$ of the intersection $\bigcap_{i\in J} E_i$ we associate the simplex $\Delta_Z:=\{w\in (\mathbb{R}_{\ge 0})^{|J|}\mid \sum_{i\in J} w_i b_i\le 1\}.$ 
This collection then defines the dual complex $\Delta_{\mathcal{X}}$ of $\mathcal{X}$. If $\mathcal{X}'$ dominates $\mathcal{X}$, we get a simplicial map from $\Delta_{\mathcal{X}'}$ to $\Delta_{\mathcal{X}}$, and one can show that there is a natural identification between $X^{\NA}$ and the projective limit of the projective system of the dual complexes and their simplicial maps. There are also natural injections $i_{\mathcal{X}}: \Delta_{\mathcal{X}}\hookrightarrow X^{\NA}$, hence one can think of $\Delta_{\mathcal{X}}$ as a subset of $X^{\NA}$, with $\mathcal{X}_0^{\com}$ being its set of vertices.

\subsection{Vertical divisors and PL functions} Recall that a divisor $D$ on a model $\mathcal{X}$ is said to be \emph{vertical} if it is supported on $|\mathcal{X}_0|$. 

Let $D$ be such a vertical divisor on $\mathcal{X}$. If $x\in X^{\di}$ we let $\mu:\mathcal{X}'\to \mathcal{X}$ be a dominating SNC model such that $x=v_{E_i}$ for an irreducible component $E_i$ of $\mathcal{X}_0'=\sum_i b_iE_i$. We write $\mu^*(D)=a_iE_i+\sum_{j\neq i} a_jE_j$ and let $f_D(x):=b_i^{-1}a_i$. This defines a function $f_D$ on $X^{\di}$ which can be seen to have a continuous extension to the whole of $X^{\NA}$, also denoted by $f_D$. Functions of this kind are called piecewise linear (PL), and and we write $f_D\in PL$.

\subsection{K\"ahler potentials and $A$-psh functions} 
A PL function $f_D$ is said to be a K\"ahler potential (with respect to $A:=\pi_X^*\alpha$), written $f_D=\phi_D\in \mathcal{H}_A$, if $A+D$  is relatively Kähler on some model $\mathcal{X}$, i.e. if $A+D+c\mathcal{X}_0$ is Kähler for large $c$.

A decreasing limit $\psi$ of K\"ahler potentials $\phi_i$ is said to be $A$-psh, and the set of $A$-psh functions is denoted by $\PSH_A$.

\subsection{Finite energy potentials} The energy $\E_A(\phi_D)$ of a K\"ahler potential $\phi_D$ is defined as $$\E_A(\phi_D):=\frac{V^{-1}}{n+1}(A+D)^{n+1},$$ while the energy $\E_A(\psi)$ of an $A$-psh function $\psi$ is defined as the infimum of the energy of all $\phi\in \mathcal{H}_A$ such that $\phi\ge \psi$. The space of finite energy potentials is defined as $$\mathcal{E}^1_A:=\{\psi\in \PSH_A: \E_A(\psi)>-\infty\}.$$ We furthermore equip $\mathcal{E}^1_A$ with the strong topology, defined as the coursest topology, finer than the topology of pointwise convergence on $X^{\di}$, such that $\E_A$ becomes continuous. We also let $\mathcal{E}^1_{A,0}:=\{\psi\in \mathcal{E}^1_A : \sup \psi=0\}$ with its subspace topology.

\subsection{The non-Archimdedean Monge-Amp\`ere operator} \label{sec:geom}

The n-A Monge-Amp\`ere measure of a K\"ahler potential $\phi_D$ is defined as $$\MA_A(\phi_D):=V^{-1}\sum_i \left((A+D)^n\cdot (b_iE_i)\right)\delta_{v_{E_i}},$$ where $\mathcal{X}_0=\sum_i b_iE_i$ is the zero divisor on a model $\mathcal{X}$ where $D$ is defined. Thus, after possibly adding a multiple of $\mathcal{X}_0$ to make $A+D$ Kähler, $\MA_A(\phi_D)$ is basically the same as $\MA(\mathcal{X},A+D)$, but thought of as a measure on $X^{\di}$ rather than on just $\mathcal{X}_0^{\com}$.

There is also a natural extension of the n-A Monge-Amp\`ere operator to $\mathcal{E}^1_A$. In that case we similarly think of the Monge-Amp\`ere measure as encoding the distribution of volume after the degeneration of $(X,\alpha)$, but insted of $X$ being decomposed into finitely many pieces, the decomposition can now also be infinitesimal. 

\subsection{Finite energy measures}
The (dual) energy $\E^{\vee}_{A}(\mu)$ of a Radon probability measure is defined as $$\E^{\vee}_{A}(\mu):=\sup\left\{E_{A}(\phi)-\int\phi\, \mathrm d\mu: \phi\in \mathcal{E}^1_{A}\right\},$$ those with $\E^{\vee}_{A}(\mu)<\infty$ giving us the space of finite energy measures $\mathcal{M}^1_{A}$. We endow $\mathcal{M}^1_{A}$ with the strong topology, defined as the coursest topology, finer than the weak topology of measures, that makes $\E_A^{\vee}$ continuous.

\subsection{A Calabi-Yau theorem for non-Archimedean finite energy potentials}

We now come to a n-A version of Theorem \ref{thm:CY}.

\begin{theorem} \label{thm:main}
The non-Archimedean Monge-Amp\`ere operator is a homeomorphism between $\mathcal{E}^1_{A,0}$ and $\mathcal{M}^1_{A}$.
\end{theorem}

As already noted in the introduction, in the algebraic setting, i.e. when $X$ is projective and $\alpha=c_1(L)$ for some ample ($\mathbb{R}$-)line bundle $L$, Theorem \ref{thm:main} was first proved by Boucksom-Jonsson \cite{BJ22}, extending earlier results of Boucksom-Jonsson-Favre \cite{BFJ15}. In the general Kähler setting, Theorem \ref{thm:main} was recently proved by Mesquita-Piccione and myself in \cite{MW25}. As an application of Theorem \ref{thm:main}, in \cite{MW25} we prove that if $(X,\alpha)$ is uniformly K-stable for models, then there is a unique cscK metric in $\alpha$. This was first proved in the algebraic case by Chi Li \cite{Li22}. See also \cite{BJ25} and \cite{DZ25} for recent related results.

\begin{rem}
In this paper we only really discuss the case when the n-A field is $\mathbb{C}$ with its trivial norm. In the algebraic setting though, much of the pluripotential theory introduced in this Section has been developed for very general n-A fields, and Calabi-Yau theorems have been proved in great generality (see e.g. \cite{BFJ15,BGJKM20,BGM22,BJ22,BJ23a} and references therein).
\end{rem}

\section{Big cohomology classes and their volumes} \label{sec:big}

\subsection{Notions of positivity for cohomology classes} Let $(X,\omega)$ be a compact K\"ahler manifold of complex dimension $n$. 

Recall that a class in $H^{1,1}(X,\mathbb{R})$ is said to be Kähler if it contains a Kähler form. The set of Kähler classes constitutes an open convex cone $\mathcal{K}:=\mathcal{K}(X)$ in $H^{1,1}(X,\mathbb{R})$ called the Kähler cone. Its closure $\overline{\mathcal{K}}$ is called the nef cone, and a class is said to be nef if it lies in the nef cone.

A class in $H^{1,1}(X,\mathbb{R})$ is said to be pseudoeffective if it contains a closed positive current. The set of pseudoeffective classes forms a closed convex cone $\mathcal{E}$ in $H^{1,1}(X,\mathbb{R})$ called the pseudoeffective cone. Its interior $\mathcal{E}^{\circ}$ is called the big cone, and a class is said to be big if it lies in the big cone. Note that a class is big if and only if it can be written as the sum of a Kähler class with a pseudoeffective class.

There are also two important notions of positivity for $(n-1,n-1)$-classes. 

A class in $H^{n-1,n-1}(X,\mathbb{R})$ is said to be pseudoeffective if it contains a closed positive current, and the set of pseudoeffective classes forms a closed convex cone $\mathcal{N}$ called the pseudoeffective cone. 

The movable cone $\mathcal{M}$ is the closed convex cone in $H^{n-1,n-1}(X,\mathbb{R})$ generated by classes of the form $\mu_*(\tilde{\beta}_1\wedge...\wedge\tilde{\beta}_{n-1})$, where $\mu:\tilde{X}\to X$ is some smooth modification and $\tilde{\beta}_i$ are Kähler classes on $\tilde{X}$. Clearly $\mathcal{M}\subseteq \mathcal{N}$.

If $X$ is projective, the Neron-Severi space $NS(X,\mathbb{R})\subseteq H^{1,1}(X,\mathbb{R})$ is defined as the subspace of $H^{1,1}(X,\mathbb{R})$ generated by divisor classes. Similarly $\mathcal{N}_1(X,\mathbb{R})\subseteq H^{n-1,n-1}(X,\mathbb{R})$ is the subspace generated by curve classes. Intersecting with $NS(X,\mathbb{R})$ or $\mathcal{N}_1(X,\mathbb{R})$ we get algebraic versions of all the positivity cones defined above. To signify the algebraic version of a cone we add the subscript $NS$. It is here important to note that a line bundle $L$ is ample/nef/pseudoeffective/big if and only if $c_1(L)$ is Kähler/nef/pseudoeffective/big.

\subsection{Duality between cones}
Note that there is a natural pairing (sometimes called the Poincaré pairing) between $H^{1,1}(X,\mathbb{R})$ and $H^{n-1,n-1}(X,\mathbb{R})$, given by $(\gamma\cdot \eta):=\int_X\gamma\wedge\eta$. 

Demailly-Păun famously proved in \cite{DP04} that $\overline{\mathcal{K}}$ and $\mathcal{N}$ are dual with respect to this pairing, i.e. that a class $\beta\in H^{1,1}(X,\mathbb{R})$ is nef if and only if $\int_X \beta\wedge \eta\ge 0$ for all $\eta\in \mathcal{N}$.

In \cite{BDPP13} Boucksom-Demailly-Păun-Peternell proved that when $X$ is projective, then $\mathcal{E}_{NS}$ and $\mathcal{M}_{NS}$ are dual with respect to the Poincaré pairing. They also conjectured that the full cones $\mathcal{E}$ and $\mathcal{M}$ should be dual for any compact Kähler manifold $X$, see \cite[Conjecture 2.3]{BDPP13}. In \cite{WN19} I proved this conjecture in the case when $X$ is projective, but the general case remains open. We will return to this issue in Section \ref{sec:TMI}, when discussing transcendental Morse inequlities. 

\subsection{$\theta$-psh functions and closed positive currents with analytic/minimal singularities} We let $\beta\in H^{1,1}(X,\mathbb{R})$ be a big class and choose a smooth form $\theta$ in $\beta$. An upper semicontinuous (usc) function $\phi: X\to[-\infty,\infty)$ is said to be $\theta$-psh if $\theta+dd^c\phi\geq0$ as a current, and the set of $\theta$-psh functions is denoted by $\PSH_{\theta}$. We say that $\phi$ has analytic singularities if locally we can write $\phi=c\log(\sum_i |f_i|^2)+g$ where the $f_i$:s are holomorphic and $g$ is bounded. If $\phi,\psi\in \PSH_{\theta}$ we say that $\phi$ is less singular than $\psi$ if $\phi\geq \psi-C$ for some constant $C$, and $\phi$ is said to have minimal singularities if $\phi$ is less singular than all $\theta$-psh functions. It is easy to see that such functions always exist. We also say that a closed positive current $T=\theta+dd^c\phi$ has analytic/minimal singularities if $\phi$ has analytic/minimal singularities.

\subsection{Lelong numbers} \label{sec:lelong} 
Let $u$ be a plurisubharmonic (psh) function in a neighbourhood of $0\in \mathbb{C}^n$. Then the Lelong number of $u$ at $0$ is defined as $$\nu_0(u):=\liminf_{z\to 0}\frac{u(z)}{\log|z|}.$$ If $u$ rather is psh in a neighbourhood of a point $x\in X$ where $X$ is a complex manifold, then $\nu_x(u):=\nu_0(u\circ g^{-1})$ where $g$ is a local holomorphic chart centered at $x$. If $T$ is a closed positive current which locally near $x$ can be written as $T=dd^c u$, then we let $\nu_x(T):=\nu_x(u)$, and if $Z\subseteq X$ is a subvariety we let $\nu_Z(T):=\inf_{x\in Z}\nu_x(T)$. 

Given a big class $\beta$ we also let $$\nu_Z(\beta):=\inf\{\nu_Z(T): T\text{ is a closed positive current in }\beta\}.$$ If $T\in \beta$ has minimal singularities then $\nu_Z(T)=\nu_Z(\beta)$. Note that Lelong numbers depend continuously on the class, see e.g. \cite[Proposition 3.6]{Bou04}.

\subsection{Non-Kähler and non-nef loci}
A closed positive current $T=\theta+dd^c\phi\in \beta$ is called a K\"ahler current if $T-\epsilon \omega\geq0$ for some $\epsilon>0$. We say that $x\in X$ lies in the K\"ahler locus of $\beta$ if there is a K\"ahler current $T\in \beta$ with analytic singularities which is smooth near $x$. The complement of the K\"ahler locus is called the non-K\"ahler locus of $\beta$ and is denoted by $E_{nK}(\beta)$.

The non-nef-locus is defined as $E_{nn}(\beta):=\{x\in X: \nu_x(\beta)>0\}$ and it is easy to see that $E_{nn}(\beta)\subseteq E_{nK}(\beta)$. It is also easy to show that if $x\in E_{nK}(\beta)\setminus E_{nn}(\beta)$ then $x$ lies in the K\"ahler locus of $\beta+\epsilon\alpha$ for any $\epsilon>0$.

\subsection{Positive products of currents} If $T_1,...,T_k$ are closed positive currents, following \cite{BEGZ10} one can form a closed positive $(k,k)$-current $\langle T_1 \wedge ... \wedge T_k\rangle$ known as the positive product. In the special case of $k=n$ and $T_i=T=\theta+dd^c\phi$ for all $i$ we have that $\langle (\theta+dd^c\phi)^n\rangle$ is a positive measure which up to a constant is equal the non-pluripolar Monge-Amp\`ere measure $\MA_{\theta}(\phi)$ (see \cite{BEGZ10}). 

\subsection{Positive intersections and volumes} \label{sec:posint} For $1\leq k\leq n$ the positive (or movable) intersection class $\langle \beta^k \rangle\in H^{k,k}(X,\mathbb{R})$ is defined as $$\langle \beta^k \rangle:=\{\langle T^k\rangle\},$$ where $T$ is any closed positive current in $\beta$ with minimal singularities. Equivalently, if $\gamma\in H^{n-k,n-k}(X,\mathbb{R})$ is semipositive (i.e. contains a semipositive form), then $\langle \beta^k \rangle\cdot \gamma$ is the supremum of all numbers $(\beta')^k\cdot \mu^*\gamma$ where $\mu: X'\to X$ is a modification and $\beta'$ is a K\"ahler class on $X'$ such that $\beta'\leq\mu^*\beta$ (see e.g. \cite{Bou02}). In the special case of $k=n$ we get that $\langle\beta^n\rangle$ is a positive number, also known as the volume of $\beta$, written $\textrm{vol}(\beta)$. If $\beta$ is K\"ahler then clearly $\langle\beta^k\rangle=\beta^k$.

If $L$ is a big line bundle, then we have that $$\langle c_1(L)^n\rangle=\textrm{vol}(L):=\limsup_{k\to\infty}\frac{h^0(X,L^k)}{k^n/n!}.$$

\subsection{Restricted volumes} \label{sec:restricted}
Let $Y$ be a smooth hypersurface on $X$ which is not contained in $E_{nK}(\beta)$. The restricted volume of $\beta$ along $Y$ is then defined as $$\langle \beta^{n-1}\rangle_{X|Y}:=\int_Y\langle (T_{|Y})^{n-1}\rangle,$$ where $T$ is any closed positive current in $\beta$ with minimal singularities. Equivalently it can be defined as the supremum of all numbers $(\beta')^{n-1}\cdot \tilde{Y}$ where $\tilde{Y}$ is the strict transform of $Y$ under a modification $\mu: X'\to X$ and $\beta'$ is a K\"ahler class on $X'$ such that $\mu^*\beta-\beta$ is the class of an effective divisor $D$ whose support is contained in $\mu^{-1}(E_{nK}(\beta))$ (see \cite[Theorem 5.3]{CT22}). 
In the case when $Y$ is contained in $E_{nK}(\beta)$ but not in $E_{nn}(\beta)$ we let $\langle \beta^{n-1}\rangle_{X|Y}:=\lim_{\epsilon \to 0+}\langle (\beta+\epsilon\{\eta\})^{n-1}\rangle_{X|Y}$, while if $Y\subseteq E_{nn}(\beta)$ we let $\langle \beta^{n-1}\rangle_{X|Y}:=0$. We say that $Y$ is $\beta$-good if it either intersects the K\"ahler locus or lies in the non-nef locus of $\beta$. For a generic big class $\beta$, $E$ will be $\beta$-good.

It follows easily from the definitions that 
\begin{equation} \label{eq:monovol}
0\le\langle \beta^{n-1}\rangle_{X|Y}\le \langle \beta^{n-1}\rangle \cdot Y\le \langle (\beta_{|Y})^{n-1}\rangle.
\end{equation}

\subsection{The transcendental Morse inequality and orthogonality} \label{sec:TMI}

In \cite{BDPP13} Boucksom-Demailly-Păun-Peternell conjectured that for any two nef classes $\alpha,\beta\in \overline{\mathcal{K}}$, the following inequality holds:
\begin{equation} \label{eq:TMI}
\langle(\alpha-\beta)^n\rangle\ge \alpha^n-n\alpha^n\cdot \beta.
\end{equation}
This is known as the transcendental Morse inequality. 

In \cite{BDPP13} it was furthermore shown that (\ref{eq:TMI}) would imply the conjectured duality between the big cone $\mathcal{E}\subseteq H^{1,1}(X,\mathbb{R})$ and the movable cone $\mathcal{M}\subseteq H^{n-1,n-1}(X,\mathbb{R})$.

When $X$ is projective, it is not hard to establish the inequality (\ref{eq:TMI}) for $\alpha,\beta\in \overline{\mathcal{K}}_{NS}$, and this is indeed how the duality between $\mathcal{E}_{NS}$ and $\mathcal{M}_{NS}$ was proved in \cite{BDPP13}.

In \cite{WN19} I proved (\ref{eq:TMI}) in the case when $X$ is projective (but obviously not assuming $\alpha$ and $ \beta$ to be algebraic). For general compact Kähler manifolds though, the transcendental Morse inequality is still a conjecture.

\subsection{Differentiability of the volume and orthogonality} \label{sec:diff}

For $X$ projective, Boucksom-Favre-Jonsson \cite{BFJ09} and Lazarsfeld-Mustaţă \cite{LM09} independently proved that the volume function $\beta\mapsto \langle\beta^n\rangle$ is continuously differentiable on the algebraic big cone $\mathcal{E}^{\circ}_{NS}$. In particular, Boucksom-Favre-Jonsson showed that if $\gamma\in NS(X,\mathbb{R})$, then
\begin{equation} \label{eq:divdiff1}
\frac{d}{dt}_{|t=0}\langle(\beta+t\gamma)^n\rangle=n\langle\beta^{n-1}\rangle \cdot \gamma,
\end{equation}
while Lazarsfeld-Mustaţă showed that if $E$ is a prime divisor not included in the non-Kähler locus of $\beta$ (which when $\beta=c_1(L)$ is the same as the non-ample locus of $L$), then 
\begin{equation} \label{eq:divdiff2}
\frac{d}{dt}_{|t=0}\langle(\beta+tE)^n\rangle=n\langle\beta^{n-1}\rangle_{X|E}.
\end{equation}
As explained in the appendix of \cite{WN19} written by Boucksom, (\ref{eq:divdiff1}) is a consequence of the Morse inequality (\ref{eq:TMI}), and since for $X$ projective, (\ref{eq:TMI}) was established in \cite{WN19} for general nef classes, it similarly follows that when $X$ is projective, the volume function is continuously differentiable on the whole of $\mathcal{E}^{\circ}$.

In the appendix of \cite{WN19} it is also explained that (\ref{eq:divdiff1}) is equivalent to the orthogonality property, established in \cite{BDPP13}, which says that for any $\beta\in \mathcal{E}^{\circ}_{NS}$, 
\begin{equation} \label{eq:orth}
\langle \beta^n\rangle=\langle \beta^{n-1}\rangle\cdot \beta.
\end{equation}

\begin{rem}
It is well known that in the algebraic case, orthogonality, or if you will the differentiability of the volume function, is tightly linked with the solvability of the n-A Monge-Amp\`ere equation, see e.g. \cite[Appendix A]{BFJ15}. This is also true when considering much more general n-A fields. What matters then is the differentiability of the associated n-A volume function, see e.g. \cite{BGJKM20,BGM22}. 
\end{rem}

\subsection{Differentiability of the volume in divisorial directions}

In \cite{WN24} I further developed the techniques of \cite{WN19} and managed to prove the following Kähler version of (\ref{eq:divdiff2}).

\begin{thm} \label{thm:diffvol}
If $Y$ is a smooth hypersurface which is $\beta$-good, then we have that 
\begin{equation} \label{eq:diff}
\frac{d}{dt}_{|t=0}\langle(\beta+tY)^n\rangle=n\langle\beta^{n-1}\rangle_{X|Y}.
\end{equation}
\end{thm}

Recently it was shown by Vu \cite{Vu23}  that just as in the algebraic case, one can remove the assumptions on $Y$ to be smooth and $\beta$-good.

An important aspect of the identity (\ref{eq:diff}) is that it shows that the restricted volume along $Y$ only depends on the cohomology class of $Y$, since this is true for the left hand side. In particular, we get the following useful corollary (see \cite[Corollary A]{WN24}):

\begin{cor} \label{cor:WN}
If $Y_1,...,Y_m$ and $Z_1,..., Z_l$ are smooth and $\beta$-good hypersurfaces such that $$\sum_ia_i\{Y_i\}=\sum_jb_j\{Z_j\},$$ then we have that $$\sum_ia_i\langle\beta^{n-1}\rangle_{X|Y_i}=\sum_jb_j\langle\beta^{n-1}\rangle_{X|Z_j}.$$
\end{cor} 

We will see that this result will play a key role in the proofs of Theorems \ref{thm:A} and \ref{thm:main}.

\section{A direct proof of Theorem \ref{thm:A}} \label{sec:direct}

We will start by showing that the Monge-Amp\`ere measure of a big test configuration is a probability measure.

\begin{proof}[Proof of Theorem \ref{thm:prob}]
Let $(\mathcal{X},A+D)$ be a big test configuration and write $\mathcal{X}_0=\sum_ib_iE_i$. Recall that $$\MA(\mathcal{X},A+D):=V^{-1}\sum_i b_i\langle(A+D)^n\rangle_{\mathcal{X}|E_i}\delta_{E_i},$$ hence we need to show that 
\begin{equation} \label{eq:voleq}
\sum_i b_i\langle(A+D)^n\rangle_{\mathcal{X}|E_i}=V.
\end{equation}
Since the restricted volumes are unaffected by adding a positive multiple of $\mathcal{X}_0$ to $D$, we can without loss of generality assume that $D\ge \mathcal{X}_0$. As a consequence, we have that $E_{nK}(A+D)\subseteq \mathcal{X}_0$, and thus $X_1$ is $(A+D)$-good.

Let us now assume that the $E_i$:s are all $(A+D)$-good. Since $$\{X_1\}=\{\mathcal{X}_0\}=\sum_i b_i\{E_i\},$$ it then follows from Corollary \ref{cor:WN} that $$\sum_i b_i\langle(A+D)^n\rangle_{\mathcal{X}|E_i}=\langle(A+D)^n\rangle_{\mathcal{X}|X_1}.$$ Since the support of $D$ does not contain $X_1$, it follows by monotonicity that $$\langle(A+D)^n\rangle_{\mathcal{X}|X_1}\ge \langle A^n\rangle_{\mathcal{X}|X_1}=A^n\cdot X_1=V.$$ On the other hand, by (\ref{eq:monovol}) we also have that $$\langle(A+D)^n\rangle_{\mathcal{X}|X_1}\le \langle(A+D)_{|X_1}^n\rangle=\langle A_{|X_1}^n\rangle=V,$$ which thus shows (\ref{eq:voleq}) as claimed. 

Assume now that some $E_i$:s are not $(A+D)$-good. Let $D'\ge \mathcal{X}_0$ be such that $A+D'$ is Kähler, and let $$B_{\epsilon}:=A+D+\epsilon(A+D').$$ Then it is easy to see that all the $E_i$:s are $B_{\epsilon}$-good for $\epsilon>0$ small enough. By the argument above we get that $$\sum_i b_i\langle B_{\epsilon}^n\rangle_{\mathcal{X}|E_i}=(1+\epsilon)^nV,$$ and we then get (\ref{eq:voleq}) by letting $\epsilon$ go to zero. 
\end{proof}

Let $\mathcal{E}^{\circ}_{\mathcal{X},A}$ denote the set of vertical divisors $D$ on a given model $\mathcal{X}$ such that $A+D$ is big.

\begin{prop} \label{prop:cont}
The functions $g_i(D):= \langle(A+D)^n\rangle_{\mathcal{X}|E_i}$ are all continuous on $\mathcal{E}^{\circ}_{\mathcal{X},A}$.
\end{prop}

\begin{proof}
Let $D\in \mathcal{E}^{\circ}_{\mathcal{X},A}$. As shown in \cite[Corollary 5.5]{CT22}, $g_i$ is continuous at the point $D$ if $E_i$ is not contained in $E_{nK}(A+D)$. If $E_i$ is contained in $E_{nn}(A+D)$, then $g_i$ is zero in a neighbourhood of $D$, and is in particular continuous at the point $D$. If $E_i$ is not $(A+D)$-good, then using monotonicity, it is easy to see that $g_i$ is continuous at $D$ if and only if $g_i(D)=0$. 

Thus let us assume that $E_i$ is not $(A+D)$-good and let also $B_{\epsilon}$ be as in the proof of Theorem \ref{thm:prob}. For $\epsilon<0$ we then have that $E_i\subseteq E_{nn}(B_{\epsilon})$ and thus $\langle B_{\epsilon}^n\rangle_{\mathcal{X}|E_i}=0$. By monotonicity we also have that for $j\neq i$, $$\langle B_{\epsilon}^n\rangle_{\mathcal{X}|E_j}\leq g_j(D).$$ Combining this with Theorem \ref{thm:prob} we get that for all $-1<\epsilon<0$, $$(1-\epsilon)^nV=\sum_j b_j\langle B_{\epsilon}^n\rangle_{\mathcal{X}|E_j}\leq \sum_{j\neq i} g_j(D)=V-g_i(D),$$
which clearly implies that $g_i(D)=0$.
\end{proof}

For each $E_i$ we have the pseudoeffective cone $\mathcal{E}_i:=\mathcal{E}(E_i)\subseteq H^{1,1}(E_i,\mathbb{R})$. Note that $\mathcal{E}_i$ does not contain any lines through the origin, because if there were such, there would have existed a non-zero class $\gamma\in H^{1,1}(E_i,\mathbb{R})$ such that both $\gamma$ and $-\gamma$ are pseudoeffective. Thus there would have been a non-trivial closed positive current of the form $dd^c\psi$, which is impossible because $E_i$ is compact. Since thus $\mathcal{E}_i$ is a closed convex cone that contains no lines, one can find a half-space $H_i\subseteq H^{1,1}(E_i,\mathbb{R})$ such that $H_i\cap \mathcal{E}_i=\{0\}$. We now pick such half-spaces $H_i$ for each $i$.

\begin{lemma} \label{lem:D_i}
For each $i$ one can find a $D_i\in \mathcal{E}^{\circ}_{\mathcal{X},A}$ such that for all $k\neq i$, $$(A+D_i)_{|E_j}\in H_k^{\circ}.$$
\end{lemma}

\begin{proof}
Let $\Gamma$ be the intersection graph of $\mathcal{X}_0$, meaning that $\Gamma$ has vertices $1,\dots,N$ and an edge between $i$ and $j$ if and only if $E_i$ and $E_j$ intersect. Let $l(i,j)$ denote length between $i$ and $j$ in $\Gamma$. 

For $1\le i\le N$ and $t\ge 0$ we let $$D_{i,t}:=\sum_j(1+t^N-t^{N-l(i,j)})b_jE_j\in \mathcal{E}^{\circ}_{\mathcal{X},A}.$$

Let $J_k:=\{j: l(j,k)=1\}$, and if $k\neq i$ we let $$J_{k,i}:=\{j: l(j,k)=1, l(i,k)-l(i,j)=1\}.$$ Note that $J_{k,i}$ is non-empty. If $j\in J_k$ we let $E_{jk}:=E_j\cap E_k$, and from the fact that ${\mathcal{X}_0}_{|E_k}={X_1}_{|E_k}=0$ we get that $$(b_kE_k)_{|E_k}=-\sum_{j\in J_k}b_jE_{jk}.$$

Combining these observations we get that 
\begin{eqnarray*}
{D_{i,t}}_{|E_k}=\sum_{j\in J_k}(t^{N-l(i,k)}-t^{N-l(i,j)})b_jE_{jk}=\\=-t^{N+1-l(i,k)}\sum_{j\in J_{k,i}}b_jE_{jk}+O(t^{N-l(i,k)}).
\end{eqnarray*}

Note that $\gamma_{k,i}:=\sum_{j\in J_{k,i}}b_jE_{jk}\in \mathcal{E}_k\setminus \{0\}$, and hence $-\gamma_{k,i}\in H_k^{\circ}$.

Now we observe that $t^{l(i,k)-N-1}(A+D_{i,t})_{|E_k}$ converges to $-\gamma_{k,i}$ as $t\to \infty$, and thus for $t_0$ large enough and $D_i:=D_{i,t_0}$, we get that $(A+D_i)_{|E_k}\in H_k^{\circ}$ for all $k\neq i$. 
\end{proof}

We are now ready to prove Theorem \ref{thm:A}.

\begin{proof}[Proof of Theorem \ref{thm:A}]
Let $\Sigma_N$ denote the unit simplex $$\Sigma_N:=\left\{a\in (\mathbb{R}_{\ge0})^N : \sum_i a_i=1\right\},$$ and pick $D_i$:s as in Lemma \ref{lem:D_i}. Given $a\in \Sigma_N$ we then let $D_a:=\sum_i a_iD_i\in \mathcal{E}^{\circ}_{\mathcal{X},A}$. We also let $f: \Sigma_N\to \Sigma_N$ be defined by $$f(a):=(V^{-1}b_1\langle(A+D_a)^n\rangle_{|E_1},\dots,V^{-1}b_N\langle(A+D_a)^n\rangle_{|E_N}).$$ By Proposition \ref{prop:cont}, $f$ is continuous. 

We now claim that $f(a)_j=0$ whenever $a_j=0$. To see this, we note that $a_j=0$ means that $D_{a}=\sum_{i\neq j}a_iD_i$, and we thus need to show that $$\langle(A+\sum_{i\neq j}a_iD_i)^n\rangle_{\mathcal{X}|E_j}=0.$$ Using the fact that $H_j^{\circ}$ is convex we get that $$(A+\sum_{i\neq j}a_iD_i)_{|E_j}=-\sum_{i \neq j}a_i\gamma_{j,i}\in H_j^{\circ}.$$ It follows that $(A+\sum_{i\neq j}a_iD_i)_{|E_j}$ is not pseudoeffective, and hence by (\ref{eq:voleq}) $$\langle(A+\sum_{i\neq j}a_iD_i)^n\rangle_{\mathcal{X}|E_j}\le \langle((A+\sum_{i\neq j}a_iD_i)_{|E_j})^n\rangle=0.$$

We clearly have that $f_t(a):=(1-t)f(a)+ta$, $t\in[0,1]$, is a homotopy of $f: \Sigma_N\to \Sigma_N$ to the identity map. Since $f(a)_j=0$ whenever $a_j=0$, each $f_t$ maps $\partial \Sigma_N$ to $\partial \Sigma_N$, thus showing that $f_{|\partial \Sigma_N}:\partial \Sigma_N\to \partial \Sigma_N$ is homotopic to the identity. If there now was a point $a_0\in \Sigma_N\setminus f(\Sigma_N)$, then we could let $g$ be a retraction of $\Sigma_N\setminus \{a_0\}$ to $\partial \Sigma_N$, and $g\circ f: \Sigma_N\to \partial \Sigma_N$ would then be a continuous extension of $f_{|\partial \Sigma_N}$ to $\Sigma_N$. But as $\Sigma_N$ is contractible, this would force $f_{|\partial \Sigma_N}$ to be null-homotopic, which it is not. Thus we must have that $f(\Sigma_N)=\Sigma_N$, which proves the theorem.
\end{proof}

\section{About the proof of Theorem \ref{thm:main}} \label{sec:proof}

Recall that Theorem \ref{thm:main} first was proved in the algebraic case by Boucksom-Jonsson in \cite{BJ22}. Their proof uses a variational method, and is close in spirit to the variational proof of Theorem \ref{thm:CY}, that was given in \cite{BBGZ13}. In \cite{BJ23a} Boucksom-Jonsson further developed their synthetic approach to non-Archimedean pluripotential theory, highlighting the similarities with classical pluripotential theory. The proof of the general Kähler version in \cite{MW25} follows the proofs in \cite{BFJ15,BJ22,BJ23a}, and as one would expect, many of the arguments can be easily translated to the Kähler setting.

Let $\mu\in \mathcal{M}^1_A$ be a probability measure on $X^{\NA}$ with finite energy, which we recall means that $$\E^{\vee}_A(\mu):=\sup\{\E_A(\phi)-\int\phi d\mu : \phi\in \mathcal{E}^1_A\}<\infty.$$ A maximizing sequence (with respect to $\mu$) is a sequence $\phi_i\in \mathcal{E}^1_A$ such that $\E_A(\phi_i)-\int\phi_i d\mu$ converges to $\E^{\vee}_A(\mu)$. A key step in solving the associated Monge-Amp\`ere equation 
\begin{equation}\label{eq:MAE}
\MA_A(\phi)=\mu,
\end{equation}
is to show that if $\phi_i$ is a maximizing sequence, then the Monge-Amp\`ere measures $\MA_A(\phi_i)$ converge strongly to $\mu$. After that one shows that some subsequence of $\phi_i$ will converge strongly to some $\phi\in \mathcal{E}^1_A$, and that this $\phi$ then solves (\ref{eq:MAE}).

As explained in \cite[Theorem 2.22]{BJ23a}, that the Monge-Ampère measures of a maximizing sequence converges strongly to $\mu$ is a consequence of the so-called orthogonality property.

\begin{df} \label{def:orth}
We say that $X^{\NA}$ has the orthogonality property (with respect to $A$) if for any continuous function $f$ on $X^{\NA}$ it holds that 
\begin{equation} \label{eq:OP}
\int(f-P_A(f))\MA_A(P_A(f))=0.
\end{equation}   
\end{df}

Here $P_A(f)$ denotes the $A$-psh envelope of $f$, defined as $$P_A(f):=\sup\{\phi\in \PSH_A: \phi\le f\},$$ and by the continuity of envelopes property we know that $P_A(f)$ is continuous and $A$-psh (see \cite[Theorem 1.4.1]{MW25}).

As shown in \cite[Appendix A]{BFJ15}, in the algebraic case, (\ref{eq:OP}) follows from the orthogonality property for big line bundles (\ref{eq:orth}), which we discussed in Section \ref{sec:diff}.

\begin{thm} \label{thm:orth}
If $(X,\alpha)$ is a compact Kähler manifold with a Kähler class $\alpha$, then $X^{\NA}$ has the orthogonality property with respect to $A:=\pi^*_X\alpha$.
\end{thm}

To prove Theorem \ref{thm:orth} we will follow the arguments in \cite{MW25}, and we will see how the results of Section \ref{sec:direct} come to use.  

Let $D$ be a vertical divisor defined on a model $\mathcal{X}$. Without loss of generality we can assume that $A+D$ is big.

Key to proving Theorem \ref{thm:orth} is the following transcendental version of \cite[Theorem 1.1 (i)]{Li23a}, which says that the Monge-Amp\`ere measure of $P_A(f_D)$ is basically the same as that of the big test configuration $(\mathcal{X},A+D)$, only extended by zero from $\mathcal{X}_0^{\com}$ to $X^{\NA}$.

\begin{thm} \label{thm:MAform}
For $\mathcal{X}$ and $D$ as above, we have that
$$\MA_A(P_A(f_D)) =V^{-1}\sum_ib_i\langle(A+D)^n\rangle_{\mathcal{X}|E_i}\delta_{v_{E_i}}.$$
\end{thm}

\begin{proof}
The proof is similar to that in \cite{Li23a}, but using the results of Section \ref{sec:direct} rather than their algebraic analogs.

Assume that $E_i$ is not contained in $E_{nK}(A+D)$. Then for any $k\in \mathbb{N}$ there is a dominating model $\mu_k:\mathcal{X}^k\to \mathcal{X}$ whose center does not contain $E_i$, a K\"ahler class $A_k$ and an effective divisor $D_k$ supported on $E_{nK}(A+D)\subseteq \mathcal{X}^k_0$ (we here identify $A+D$ with its pullback to $\mathcal{X}^k$) such that $A+D=A_k+D_k$ and 
\begin{equation} \label{eq:intP}
(A_k)^n\cdot \tilde{E}_i\ge \langle (A+D)^n\rangle_{\mathcal{X}|E_i}-1/k.
\end{equation}
Here $\tilde{E}_i$ denotes the strict transform of $E_i$. We also note that $D_k$ is vertical and hence $\phi_k:=\phi_{D-D_k}\in \mathcal{H}_A$, and since $D_k\ge 0$ we have that $\phi_k\le f_D$.

Using the continuity of $P_A(f_D)$ and a simple monotonicity argument we can assume that $\phi_k$ converges uniformly to $P_A(f_D)$. The Monge-Amp\`ere operator is easily seen to be continuous with respect to uniform convergence, and thus $\MA_A(\phi_k)$ will converge weakly to $\MA_A(P_A(f_D))$. Using (\ref{eq:intP}) this implies that  
$$\MA_A(P_A(f_D))(\{v_{E_i}\})\ge \liminf_{k\to \infty} \MA_A(\phi_k)(\{v_{E_i}\})\ge V^{-1}\langle (A+D)^n\rangle_{\mathcal{X}|E_i}.$$

If $E_i$ is not $(A+D)$-good, then we see from the proof of Proposition \ref{prop:cont} that $\langle (A+D)^n\rangle_{\mathcal{X}|E_i}=0$. Thus we get the inequality $$\MA_A(P_A(f_D)) \geq V^{-1}\sum_ib_i\langle(A+D)^n\rangle_{\mathcal{X}|E_i}\delta_{v_{E_i}}.$$ By Theorem \ref{thm:prob} we know that the right hand side is a probability measure just as the left hand side, which means that they have to be equal.
\end{proof}

A probability measure $\mu$ on $X^{\NA}$ of the form $\mu=\sum_{i=1}^N a_i\delta_{v_{E_i}}$, where each $v_{E_i}\in X^{\di}$, is called a divisorial measure. Note that the fact that one can solve the n-A Monge-Amp\`ere equation for all divisorial measures, see \cite[Proposition 8.6]{BFJ15} and \cite[Proposition 8.3.2]{MP24}, now follows as a direct corollary of Theorem \ref{thm:A} and Theorem \ref{thm:MAform}.

We will also have use for the following interpretation of the numbers $f_D(v_{E_i})-P_A(f_D)(v_{E_i})$ in terms of Lelong numbers of the class $A+D$.

\begin{prop} \label{prop:lelongineq}
We have that $$f_D(v_{E_i})-P_A(f_D)(v_{E_i})= b_i^{-1}\nu_{E_i}(A+D).$$
\end{prop}
\begin{proof}
By Demailly regularization, for every $\epsilon>0$ we can find a K\"ahler current $T\in A+D$ with analytic singularities such that $\nu_{E_i}(T)\le \nu_{E_i}(A+D)-\epsilon$. Thus on some model $\mathcal{X}$, $T$ will have divisorial singularities, and since $E_{nK}(A+D)\subseteq |\mathcal{X}_0|$ we can assume that the singularities are supported on $|\mathcal{X}_0|$. Hence we can write $T=\Omega+[D']$ where $\Omega$ is K\"ahler and $D'$ is an effective vertical divisor. It follows that $\phi:=\phi_{D-D'}\in \mathcal{H}_A$ and $\phi\le f_D$, which implies that $\phi\le P_A(f_D)$. We finally get that 
\begin{eqnarray*}
f_D(v_{E_i})-P_A(f_D)(v_{E_i})\le f_D(v_{E_i})-\phi(v_{E_i})=\\= b_i^{-1}\nu_{E_i}([D'])=b_i^{-1}\nu_{E_i}(T)\le \nu_{E_i}(A+D)-\epsilon.
\end{eqnarray*}
Since $\epsilon>0$ was arbitrary this implies that $$f_D(v_{E_i})-P_A(f_D)(v_{E_i})\le b_i^{-1}\nu_{E_i}(A+D).$$

The argument for the other inequality is not very difficult, but as it is not actually needed for proof of Theorem \ref{thm:orth}, we will omit the details. 
\end{proof}

We are now ready to prove Theorem~\ref{thm:orth}.

\begin{proof}[Proof of Theorem~\ref{thm:orth}]
We first assume that $f=f_D$ is PL, and without loss of generality we can assume that $A+D$ is big. By Theorem \ref{thm:MAform} $\MA_A(P_A(f_D))$ is supported on $\mathcal{X}_0^{\com}$. Note that $\langle(A+D_{\epsilon})^n\rangle_{\mathcal{X}|E_i}=0$ if $\nu_{E_i}(A+D)>0$, so by Proposition \ref{prop:lelongineq}, $f_D-P_A(f_D)=0$ on the support of $\MA_A(P_A(f_D))$, proving the claim in the PL case.

For a general continuous function $f$ we choose a sequence of PL functions $f_i$ converging uniformly to $f$. It is then easy to see that $P_A(f_i)$ also will converge uniformly to $P_A(f)$, and thus $\MA_A(P_A(f_i))$ will converge weakly to $\MA_A(P_A(f))$. This then shows that $$\int (f-P_A(f))\MA_A(P_A(f))=\lim_i\int (f_i-P_A(f_i))\MA_A(P_A(f_i))=0,$$ which concludes the proof.
\end{proof}

\end{document}